\documentclass[12pt,a4paper]{amsart}

\usepackage{amsmath,eucal}
\usepackage{amssymb}

\newcounter{prp}
\newcounter{thm}
\newcounter{cor}
\newcounter{lem}
\newcounter{exl}

\newtheorem{lem}[lem]{Lemma}
\newtheorem{prp}[prp]{Proposition}
\newtheorem{thm}[thm]{Theorem}

\newtheorem{clm}{Claim}

\newcommand{\tg}{\tilde{g}}

\title{ Collapse of unit horizontal bundles equipped with a  metric of Cheeger-Gromoll type}
\author{Wojciech  Koz\l owski \and Szymon M. Walczak}

\begin{document}

\begin{abstract}
We study unit horizontal bundles associated with Riemannian submersions.
First we investigate metric properties of  an arbitrary unit horizontal bundle equipped with a
Riemannian metric of the Cheeger-Gromoll type.
Next we examine it from the Gromov-Hausdorff convergence theory point of view,
and we state a collapse theorem
for unit horizontal bundles
associated with a sequence of warped Riemannian submersions.
\end{abstract}

\maketitle

\section{Introduction and preliminaries}

\subsection{Introduction}

Recently in \cite{BLW}, M. Benyounes, E. Loubeau and C. M. Wood introduced a new class of natural  metrics of Cheeger-Gromoll type on
the vector bundle over a Riemannian manifold.
These metrics, ${\bf h}_{p,q}$, $p,q\in \mathbb{R}$, $q\ge 0$, called {\em $(p,q)$-metrics},
generalize Sasaki  metric  \cite{Kw}  and Cheeger-Gromoll metric \cite{GK} on $TM$.

 Although $(p,q)$-metrics have been discovered together with some new harmonics maps,  the geometry of
$(p,q)$-geometry of the tangent bundle is of the independent interest \cite{BLW1}.

In the present paper we combine a technique of $(p,q)$-metrics and Riemannian submersions with the Gromov-Hausdorff distance theory
(GH-distance theory).

First, we investigate the unit horizontal bundle $\tilde E^1$, associated with a Riemannian submersion
$P:\tilde M \to M$.  The total space of $\tilde E^1$  consists of the all unit vectors in $T\tilde M$ which are orthogonal to the fibres of $P$.

The differential $\tilde P= P_\ast$ maps $\tilde E^1$ into $SM$ - the unit sphere bundle over $M$.

Let us equip  $\tilde E^1$ and $SM$
with the  $(p,q)$-metric $\tilde {\bf h}$ and ${\bf h}$.
We ask when $\tilde P: (\tilde E^1, \tilde {\bf h}) \to ( SM, {\bf h})$ is a Riemannian submersion.
We prove (Proposition \ref{Prp_E1toSM}) that $\tilde P: \tilde E^1 \to  SM$ is a Riemannian submersion iff the horizontal distribution of $P$ is integrable.
This assertion seems to be of independent interest.

Next, we combine  the result from Proposition \ref{Prp_E1toSM}
with Theorem 2 from \cite{W},  and we obtain Collapse Theorem (Theorem \ref{Thm_Col_Thm})
for the unit horizontal bundle.  We prove that the sequence of unit horizontal bundles
$( E^1_n)_{n\in \mathbb{N}}$ associated with the sequence of
warped Riemannian submersions
$(P_n : \tilde M_n \to M)_{n \in \mathbb{N}}$ converges (in GH-topology)
to the unit sphere  bundle $SM$ iff $(\tilde M_n)_{ n\in \mathbb{N}}$ converges to $M$.

In fact, Theorem \ref{Thm_Col_Thm} asserts that this natural construction is continuous in the GH-topology. Some other examples showing the continuity
of a natural constructions  
can be found in  \cite{L} by H. Li, where the author proves continuity of  {\em $\theta$-deformations}, and
in the recent paper of P. G. Walczak \cite{WP} where  the {\em continuity of spaces of probability measures associated with a Riemannian manifolds}
is studied.

\subsection{Riemannian submersions}

\setcounter{lem}{0}
\setcounter{thm}{0}
\setcounter{prp}{0}
\setcounter{equation}{0}
\setcounter{cor}{0}
\setcounter{clm}{0}

We briefly review basic facts of Riemannian submersions. For more details we refer to \cite[Ch. 9]{Be}.

If $ P: (\tilde M, \tg)  \to (M, g)$ is a Riemannian submersion then its tangent
bundle
$T \tilde M$
splits as a direct orthogonal sum
$T\tilde M = \mathcal H \oplus \mathcal V$,
where $\mathcal V =   \ker P_\ast$ is the vertical subbundle and $\mathcal H = \mathcal V^\perp$ is the horizontal subbundle. If $W \in T\tilde M$
then $W= \mathcal H W + \mathcal V W$ denotes the corresponding orthogonal splitting.

To indicate a submersion we work with, we often write
$\mathcal H^P$ and $\mathcal V^P$ instead of $\mathcal H$ and
$\mathcal V$.

A vector field $\tilde X$
(resp. $\tilde U$) is horizontal (resp. vertical) if $\tilde X \in \Gamma(\tilde M, \mathcal H^P)$
(resp. $\tilde U\in \Gamma (\tilde M, \mathcal V^P $)).
The vector field $\tilde X$ is called {\em basic} if $\tilde X$ is horizontal and there exists a vector filed $X$ on $M$ such that
$P_\ast \tilde X = X$.  There is one-to-one correspondence $X\to \tilde X$ between vector fields on $M$ and basic vector fields on $\tilde M$. If $X$ is a vector field on $M$ then the corresponding basic vector field $\tilde X$ is called the {\em horizontal lift} of $X$.
If the vector fields $\tilde X$, $\tilde Y$ are basic and the vector field $\tilde U$ is vertical then $\tg(\tilde X, \tilde Y)$ is constant along the fibres of $P$, and $[\tilde X, \tilde U]$ is also  vertical. Moreover, the horizontal part of the Lie bracket $[\tilde X, \tilde Y]$ coincide with the horizontal lift $\widetilde{[X,Y]}$.

Let $\tilde \nabla$ and $\nabla$ be the Levi-Civita connections of $\tilde g$ and $g$, respectively. The connection $\tilde D$ in $\mathcal H^P$ is
simply the horizontal projection  $\mathcal H \tilde \nabla$.
Since $\tilde\nabla$ is  Riemannian, so $\tilde D$ is.

\begin{lem}\label{Lemma_i_ii}
Suppose that $\tilde X$, $\tilde Y$ are basic and $\tilde U$ is a vertical vector field.
Put $X=P_\ast \tilde X$ and $Y=P_\ast \tilde Y$. Then
\begin{itemize}
\item[(i)] $P_\ast \tilde D_{\tilde X} {\tilde Y} = \nabla_{ X}  Y$, and $\tilde D_{\tilde X} {\tilde Y} = \widetilde{\nabla_{ X}  Y}$
\item[(ii)] $\tg ( \tilde D_{\tilde U} {\tilde X}, \tilde Y ) = -\frac{1}{2} \tg ( \tilde U, [ \tilde X, \tilde Y ] ).$
\end{itemize}
\end{lem}

\begin{lem}\label{net_i_ii}
Let us suppose that $\tilde M$ is compact.
\begin{itemize}
\item[(i)]
Let $\{x_1,\dots, x_k\}$ be an $\varepsilon$-net in $M$, and let for every $i=1,\dots,k$, $\{\tilde x_{i1},\dots, \tilde x_{il(i)}\}$
be an $\varepsilon$-net in the fibre $\tilde M_{x_i} = P^{-1}(x_i)$. Then the set
$\{ \tilde x_{ij} : i=1,
\dots, k, j=1,\dots l(i)\}$ is a
 $(2\varepsilon)$-net in $\tilde M$.
 \item[(ii)] Let $\{ \tilde x_1, \dots, \tilde x_s\}$ be  an $\varepsilon$-net in $\tilde M$. Then
 the set $P\big(\{ x_1,\dots, x_s\}\big)$ is an $\varepsilon$-net in $M$ (notice that $P(\tilde x_i)$ may be equal to $P(\tilde x_j)$ even if $i\ne j$).
\end{itemize}
\end{lem}

We omit elementary proofs. \\

\subsection{(p,q)-metrics}

Let  $(M^n,g)$
be an arbitrary  Riemannian manifold
  with the Levi-Civita connection $\nabla$.

Let $ E$   be a vector bundle  $\pi:\mathcal E\to M$ equipped with
a fibre metric $h$ and
 a Riemannian connection $D$.

 We define  the {\em connection map }  $K=K^D: T\mathcal E\to \mathcal E$  related to $D$ as follows:
$K$ is a smooth map inducing for every  $\zeta \in \mathcal E$ a
$\mathbb{R}$-linear map $ T_\zeta \mathcal E_x\to \mathcal E_x$, $x=\pi\zeta$ and
determined by the condition:
${K}(\xi_\ast v) = D_v \xi$, $v \in TM$, $\xi\in \Gamma (M,\mathcal E)$.
Notice that by the definition follows that
$K|T_\zeta \mathcal E_x$ is the canonical isomorphism $ T_\zeta \mathcal E_x\to \mathcal E_x$.
For more details on the connection map we refer to \cite{Rz}  and \cite{Kw}.

Following by \cite{BLW},
one can equip $T\mathcal E$ with the $(p,q)$-metric  ${\bf h}_{p,q}$,  $p,q\in \mathbb{R}$, $q\ge 0$,
that is, a Riemannian metric defined by:
\[  {\bf h}_{p,q}( A, B)   = g( \pi_\ast  A, \pi_\ast B ) +\frac{1}{(1+|\zeta|^2)^p}
(h( K A,  K B ) + q h( KA,\zeta) h(KB, \zeta)),\]
where $ A, B \in T_\zeta \mathcal E$, and $|\zeta|^2 = h(\zeta,\zeta)$.

Suppose  $E=(TM, g, \nabla)$.
If $p=q=0$ then ${\bf h}_{p,q}$ coincide with the Sasaki metric \cite{Kw}.
If  $p=q=1$ then ${\bf h}_{p,q}$ becomes the Cheeger-Gromoll metric \cite{GK}.

  The projection $\pi: (\mathcal E,{\bf h}_{p,q} )\to (M,g)$  is  a Riemannian submersion
such that $\mathcal H^{\pi} = \ker K^D$ and $\mathcal V^\pi = \ker \pi_\ast$.
Consequently, for every $w\in T_x M$ and $\zeta\in E_x$ there exists a unique horizontal (resp. vertical) vector
$w^h\in \mathcal{H}_\zeta$ (resp. $w^v\in \mathcal{V}_\zeta$), i.e.,
 $\pi_\ast w^h =w$ and  $Kw^h=0$ (resp. $\pi_\ast w^v =0$ and  $Kw^v=w$).

Let $\pi^1: \mathcal E^1\to M$ be the {\em unit  bundle induced from} $\pi: \mathcal E\to M$
, i.e., $\mathcal E^1=\{\xi\in \mathcal E: |\xi|=1\}$  and $\pi^1 = \pi|\mathcal E^1$. If $\mathcal E=TM$ then
$\pi^1: \mathcal E^1\to M$ is simply the {\em unit sphere bundle} $SM$.
If $Q:M\to N$  is a Riemannian submersion
and $\mathcal E = \mathcal H^Q$ then
$\pi^1: \mathcal E^1\to M$ is called a {\em unit horizontal bundle}.

We define a $(p,q)$-metric on $\mathcal E^1$ by the restriction of ${\bf h}_{p,q}$ to $T(\mathcal E^1)$. We denote it also by ${\bf h}_{p,q}$.
Notice that $\pi^1 = \pi | \mathcal E^1 : (\mathcal E^1, {\bf h}_{p,q}) \to (M, g)$ is a Riemannian submersion
such that, for any $\xi \in \mathcal E^1$,  $\mathcal{H}^{\pi^1}_\xi =\mathcal H^\pi_\xi$ and
$\mathcal V^{\pi^1}_\xi  = \{A\in \mathcal V^{\pi}_\xi : \langle K^D A, \xi \rangle =0\}$
(cf. \cite[\S 2, Lemma 1]{Rz}).

Since ${\bf h}_{p,q}$ depends on $E=(\pi : \mathcal E \to M, h, D)$ it is convenient to identify $E$ with its total space  $\mathcal E$ and
write $(E,{\bf h}_{p,q})$ instead of $(\mathcal  E,{\bf h}_{p,q})$. Moreover,
the corresponding unit horizontal bundle with the induced $(p,q)$-metric is denoted by $(E^1, {\bf h}_{p,q})$.

\section{Results}
\subsection{Submersion theorem}
Let $P:(\tilde M, \tg)  \to (M, g)$ be a Riemannian submersion, $b=\dim M$, $a+b=\dim \tilde M$ and  $T \tilde M = \mathcal{H}^P \oplus \mathcal{V}^P$
be the corresponding orthogonal splitting of $T \tilde M$.
Let $\tilde \nabla$ and $\nabla$ denote the Levi-Civita connections of $\tilde g$ and $g$,
respectively.

Let $\tilde E$ denote the horizontal bundle $\mathcal H^P\to \tilde M$ equipped with the fibre metric $\tilde g$ and the Riemannian
connection
$\tilde D = \mathcal H \tilde \nabla$.

Let $\tilde \tau $ be the parallel transport in $\tilde E$, $K^{\tilde D}$ and
$K^\nabla$ be the  connection maps corresponding to $\tilde D$ and $\nabla$.
Moreover, let $\tilde \pi $ and $\pi $ denote the restrictions of the natural projections   $\tilde E \to \tilde M$ and $ T M \to  M$ to
$\tilde E^1$ and $SM$, respectively.

For $w\in T_x M$ and $\tilde x \in P^{-1}(x)$, let $\tilde w = \tilde w_{\tilde x} \in \mathcal{H}^P_{\tilde x}$
be the unique horizontal vector such that $P_\ast \tilde w = w$.
Similarly, for any $W \in T_{\tilde x} \tilde M$ and any $\xi \in \tilde E^1_{\tilde x}$ let $W^h\in \mathcal{H}^{\tilde \pi}_\xi$ be the unique
horizontal vector such that $\tilde \pi_\ast (W^h) = W$.

If $\tilde x\in \tilde M$,
 $P(\tilde x)=x$ and $\xi \in \tilde E^1_{\tilde x}$ then the horizontal lift
$\tilde w^h_\xi \in \mathcal{H}^{\tilde \pi}_\xi$ of $w\in T_x M$ to the point $\xi$
 may be constructed as follows: Let $\gamma$ be a curve in $M$ such that $\gamma(0)=x$ and $\dot \gamma (0) =w$.
Next, let $\tilde \gamma = \tilde \gamma_{\tilde x} $ be the horizontal lift of $\gamma$ such that $\tilde \gamma(0) = \tilde x$. Let
$ \tilde \gamma^h (t) = \tilde \gamma_\xi^h (t) = \tilde \tau^{\tilde \gamma}_t \xi$ be the
parallel transport of $\xi$ along $\tilde \gamma$ from $0$ to $t$.
Then $ \tilde w^h_\xi = \dot{\tilde\gamma}^h(0)$.

Suppose that $p,q \in \mathbb{R}$, $q\ge 0$ are fixed.
Let ${\bf h}= {\bf h}_{p,q}$ and ${\tilde{ \bf  h}} = {\tilde{ \bf  h}}_{p,q}$ denote $(p,q)$-metrics on $SM$ and $\tilde E^1$, respectively.

We distinguish the following pairwise orthogonal subbundles of $T(\tilde E^1)$:

\begin{eqnarray*}
H'_\xi (\tilde E^1) &=&  T_\xi \tilde E^1_{\tilde \pi (\xi )},\\
 H''_\xi (\tilde E^1) &=& \{\tilde w^h_\xi: P\tilde \pi(\xi)=x, w\in T_x M\},\\
 V_\xi (\tilde E^1) &=& \{ W^h \in \mathcal{H}^{\tilde \pi}_\xi : P_\ast \pi_\ast W^h =0\}.
\end{eqnarray*}

Let $\tilde P = P_\ast$. Since $P$ is a Riemannian submersion we see that
\[\tilde P= P_\ast: \tilde  E^1 \to SM. \]
We write
$H' = H'(\tilde E^1)$, $H''= H''(\tilde E^1)$ and $V=V(\tilde E^1)$ for simplicity.

\begin{clm}\label{clm1} $\tilde P_\ast : H'_\xi
\to T_u (S_xM)$, $u=\tilde P \xi$, is an isometry.
\end{clm}

\begin{proof}  Clearly, $\tilde P_\ast ( H'_\xi)
\subset T_u S_xM$ and
$\dim H'_\xi = \dim T_u (S_xM)$. Thus it suffices to show that $\tilde P_\ast $ preserves the length of vectors.

Let $\tilde \pi \xi = \tilde x$ and $\pi u = x$.
Take $A\in H'_\xi$. Since $A \in \mathcal V^{\tilde \pi}_\xi$ and $\tilde P_\ast A \in \mathcal V^\pi_u$,
\begin{eqnarray*}
|A|^2 &=& \frac{1}{(1+|\xi|^2)^p}\Big(|K^{\tilde D} A|^2 + q \big(\tilde g (K^{\tilde D}A,\xi)\big)^2\Big),\\
|\tilde P_\ast A|^2 &=& \frac{1}{(1+|u|^2)^p}\Big(|K^{\nabla} \tilde P_\ast A|^2 + q \big( g (K^\nabla \tilde P_\ast A,u)\big)^2\Big).
\end{eqnarray*}
Since $\tilde P:\tilde  E_{\tilde x} \to T_x M$ is a linear map, $K^\nabla \tilde P_\ast = P_\ast K^{\tilde D}$ on
$T_\xi \tilde E^1_{\tilde x}$. Thus
\begin{eqnarray*}
 |K^\nabla \tilde P_\ast A| &=& |P_\ast K^{\tilde D} A| = |K^{\tilde D} A|,\\
 g (K^\nabla \tilde P_\ast A,u) & =&  g (P_\ast K^{\tilde D} A, P_\ast \xi) = \tilde g ( K^{\tilde D} A,  \xi).
\end{eqnarray*}
Consequently, $|A|=|\tilde P_\ast A|$. This proves the assertion.
\end{proof}

\begin{clm}\label{clm2} $\tilde P_\ast : H''_{\xi} \to \mathcal{H}^\pi _u$, $u = P_\ast \xi$ is an isometry.
\end{clm}

\begin{proof} Let $x= \pi u$. By the definition of $H''$,  $\dim H''_{\xi} =  \dim T_x M$. On the other hand,
 $\dim \mathcal{H}^\pi _u = \dim T_x M$.  To prove the assertion it suffices to show that $\pi_\ast \tilde P_\ast $ preserves the length of vectors and
its image is contained in   $\mathcal{H}^\pi _u$.

Let $\tilde w^h\in H''_{\xi}$. We can suppose that $\tilde w^h = \dot{\tilde \gamma}^h(0)$ where $\gamma$ is a curve in $M$ such that
$w=\dot \gamma(0)$.
Then
\[|\tilde w^h|=|\tilde \pi_\ast \tilde w^h| = |\tilde w|=|P_\ast \tilde w|=|w|.\]
 Observe that $\tilde P \tilde \gamma^h$ is a vector field along $\gamma$. Indeed, we have
\[\pi \tilde P \tilde \gamma^h = \pi P_\ast \tilde \gamma^h = P \tilde \pi \tilde \gamma^h= P\tilde \gamma =\gamma.\]
Next we have
\[ |\pi_\ast \tilde P_\ast \tilde w^h|=|(\pi  \tilde P {\tilde\gamma}^h)\dot{}(0)| =|\dot \gamma(0)|=|w|.\]
Consequently, we have shown that $|\pi_\ast \tilde P_\ast \tilde w^h |=| \tilde w^h|$.
Thus $\pi_\ast \tilde P_\ast$ preserves the length of vectors belonging to  $H''_{\xi}$.

We have to show that the vertical part of $\tilde P_\ast \tilde w^h$ is equal to zero, or equivalently, $K^\nabla (\tilde P_\ast \tilde w^h)=0$.
Since $\tilde\gamma^h(t)\in \mathcal{H}^P$, by Lemma \ref{Lemma_i_ii}(i) we
conclude that $\nabla_{P_\ast \dot{\tilde\gamma}} P_\ast \tilde\gamma^h=P_\ast \tilde D_{\dot{\tilde\gamma}} \tilde\gamma^h.$
Since $\tilde P \tilde \gamma^h $ is a vector field along $\gamma$ we obtain
\[
K^\nabla (\tilde P_\ast \tilde w^h)=(\nabla_{P_\ast \dot{\tilde\gamma}} P_\ast \tilde\gamma^h)(0)=
(P_\ast \tilde D_{\dot{\tilde\gamma}}\tilde\gamma^h)(0)=0.
\]
Consequently, $\tilde P_\ast \tilde w^h\in \mathcal H^\pi_u$.
\end{proof}

 \begin{clm}\label{clm3} The following conditions are equivalent:
 \begin{itemize}
 \item[(i)]
$\tilde P_\ast(  V) =0$.
\item [(ii)]
For every fibre $P^{-1}(x)$ and every curve  $\eta$ in  $P^{-1}(x)$,
$\tilde P \tilde \tau^\eta =\tilde P$.
\item[(iii)]
$\mathcal H^P$ is integrable.
 \end{itemize}
\end{clm}

 \begin{proof} (i) $\Leftrightarrow$ (ii). Let $W^h \in   V_\xi$, $\xi \in \tilde E^1$, $\tilde \pi \xi = \tilde x$ and let
$\tilde x \in P^{-1}(x)$.
Then $W^h = (t \mapsto \tilde \tau^{\eta}_t
\xi)\dot{}(0)$ where
$\eta$ is a curve in the fibre $P^{-1}(x)$ such that $\eta(0)= \tilde x$ and $\dot \eta (0) =\tilde \pi_\ast W^h$.
Thus $\tilde P_\ast W^h= ( t \mapsto P_\ast    \tilde \tau^{\eta}_t  \xi)\dot{}(0)$. Consequently,
we see that $\tilde P_\ast W^h =0 $ for every $W^h \in V$ iff the curve
$t \mapsto P_\ast    \tilde \tau^{\eta}_t  \xi$ is constant for every $\xi$ and $\eta$.  The last condition is equivalent to  $P_\ast    \tilde \tau^{\eta} = P_\ast$.

(ii) $\Leftrightarrow$ (iii).
By Lemma \ref{Lemma_i_ii} (ii) $\mathcal H^P$ is integrable iff every basic vector field $\tilde X$ is parallel $(\tilde D \tilde X = 0$).
This is equivalent to the condition
$\tilde \tau^\eta \tilde X = \tilde X$, which is equivalent to (ii).
\end{proof}

Observe that
\begin{eqnarray}
\dim H' +\dim H''
+\dim V & =&  (b-1) + b + a =
\dim \tilde E^1,\label{wymiar1}\\
\dim H'  +\dim H'' & =&  \dim SM.\label{wymiar2}
\end{eqnarray}

\begin{lem}\label{Lemat_alg} Suppose $X$ and $Y$ are two finite dimensional Euclidean spaces. Suppose $Z$ is a subspace of $X$ and $X= Z\oplus Z^\perp $ is the corresponding orthogonal splitting. Given a linear operator $A: X\to Y$ such that  $A: Z \to Y$ is an isometry and $AZ^\perp = 0$. If $Z'$ is any linear subspace of $X$ such that $A:Z'\to Y$ is an isometry then $Z'=Z$.
\end{lem}

\begin{prp}\label{Prp_E1toSM} $\tilde P : (\tilde E^1 , \tilde {\bf h} )\to (SM, {\bf h})$ is a Riemannian submersion iff $ \mathcal{H}^P$ is integrable. Then
$\mathcal H^{\tilde P} = H'\oplus H''$ and $\mathcal V^{\tilde P} = V$.
\end{prp}

\begin{proof}
This assertion follows directly by Claims \ref{clm1}-\ref{clm3}, \eqref{wymiar1}-\eqref{wymiar2} and Lemma \ref{Lemat_alg}
\end{proof}

\subsection{Collapse theorem}\label{Section_Col_Thm}

\setcounter{lem}{0}
\setcounter{thm}{0}
\setcounter{prp}{0}
\setcounter{equation}{0}
\setcounter{cor}{0}
\setcounter{clm}{0}

Let $P: (\tilde M, \tg)\to (M,g)$ be a Riemannian submersion and let $f:M\to (0,+\infty)$ be a smooth function.  Put
$\tilde f = f \circ  P$.
We  modify the Riemannian metric $\tg$ to $\tg_f$
putting:
\begin{eqnarray*}
\tg_f &=& \tg
\text{ on }
\mathcal{H}^P\times \mathcal{H}^P,\\
\tg_f &=& \tilde f^2 \tg \text{ on }
\mathcal{V}^P\times T\tilde M.
\end{eqnarray*}
Then $P: (\tilde M, \tg_f)  \to (M, g)$ remains a Riemannian submersion. We  call it a {\em warped submersion}, while the function
$f$ is called {\em warping function} \cite{W}.
Denote by $M_f$ the Riemannian manifold $(M,g_f)$.

Notice that the horizontal distributions of $P: (\tilde M, \tg)  \to (M, g)$ an $P: (\tilde M, \tg_f)  \to (M, g)$ coincide, and are equal to
$\mathcal   H^P$.\\

The following {\em Collapse Theorem } is proved in \cite[Thm. 2]{W}.  In the whole section `$\tt lim$' denotes
the limit in the GH-topology. For basic facts concerning GH-topology we refer to \cite{P1} or \cite[\S 12.4]{Br}

\begin{thm}\label{Thm_Szymon}  Let $P:\tilde M \to M$ be a Riemannian submersion with $\tilde M$ compact. Suppose that
$(f_n: M \to (0,+\infty))_{n\in \mathbb{N}}$
is a uniformly bounded sequence of warping functions. Put $\tilde M_n= \tilde M_{f_n}$. We have $\lim \tilde M_{n} = M$
iff
for every $\varepsilon>0$ there exists a positive integer $N$ such that for every $n>N$ there exists an $\varepsilon$-net $A^{(n)}$ on $M$
such that $f_n |A^{(n)} < \varepsilon.$
\end{thm}

Suppose  a Riemannian submersion $P:(\tilde M, \tg)  \to (M, g)$  and  a warping function $f$ are given..
 Let $p,q \in \mathbb{R}$ and $q\ge 0$.

Let us denote by $\tilde E $ the vector bundle  $\tilde \pi: \mathcal H^{P}\to \tilde M$ with the fibre metric $\tilde g$.

Let
$\tilde E_f$  denote the vector bundle  $\tilde \pi: \mathcal H^{P}\to \tilde M$  with the fibre metric $\tilde g_f$.

Let $\tilde D$ and $\tilde D^f$ be the connections  in $\tilde E$ and $\tilde E_f$, respectively. 
The corresponding connections maps are denoted by $K$ and $K^f$, respectively.

 As in the previous section, let  ${\bf h}$ and $\tilde{\bf h}$ denote the $(p,q)$-metrics on $SM$ and $\tilde E^1$, respectively.
Moreover, let $\tilde{\bf h}_f$ denote the $(p,q)$-metric on $\tilde E^1_f$. We want to find relations between Riemannian manifolds 
$(\tilde E^1, \tilde{\bf h})$ and
$(\tilde E_f^1, \tilde{\bf h}_f)$.

Let $\tilde P = P_\ast$.

\begin{clm} $\tilde P : (\tilde E_f^1, \tilde{\bf h}_f) \to (SM, {\bf h})$
 is a Riemannian submersion iff $\tilde P : (\tilde E^1, \tilde{\bf h}) \to (SM, {\bf h}))$ is.
Then the vertical distributions of these bundles coincide and are equal to $V(\tilde E)$.
\end{clm}

\begin{proof} The first statement of the assertion follows directly from Proposition \ref{Prp_E1toSM} and the fact that the horizontal
 distribution of a warped submersion and the initial submersion coincide.

The second one follows from the fact that the vertical distribution is equal to the kernel of $\tilde P_\ast$.
\end{proof}

\begin{clm} We have $H'(\tilde E_f)= H'(\tilde E)$ and $H''(\tilde E_f)= H''(\tilde E)$. Consequently, if $\mathcal H^p$ is
integrable then the horizontal distributions of the Riemannian submersions
$\tilde P : (\tilde E_f^1, \tilde{\bf h}_f) \to (SM, {\bf h}))$ and
$\tilde P : (\tilde E^1, \tilde{\bf h}) \to (SM, {\bf h}))$ coincide.
\end{clm}
\begin{proof}
The first identity is obvious. The second follows by Lemma \ref{Lemma_i_ii} (i).
More precisely, suppose $\gamma$ is a curve in $M$ and $\tilde \gamma $ is its horizontal lift. Take any parallel vector field $\sigma$ along $\gamma$.
Let $\tilde \sigma$ be the unique horizontal vector field along
$\tilde \gamma$ such that
$P_\ast \tilde \sigma = \sigma$. Then
Lemma \ref{Lemma_i_ii} (i) implies that $\tilde \sigma$ is parallel vector field with respect to both connections. It follows that
$H''(\tilde E_f)= H''(\tilde E)$.
\end{proof}

\begin{clm}\label{CTClm3}
Suppose that $\mathcal H^P$ is integrable.
Then $\tilde P : (\tilde E^1_f, \tilde{\bf h}_f) \to (SM, {\bf h}) $  is a warped submersion whose warping function is
$ \hat f =f\circ \pi$, where $\pi: SM \to M$ is the natural projection.
\end{clm}

\begin{proof}
Let $\tilde \xi \in \tilde E_f$, $\tilde \pi \tilde \xi = \tilde x$, $P_\ast \tilde \xi =\xi$ and $P \tilde x = x$,
and $A,B\in T_{\tilde \xi} \tilde E^1_f$.

Suppose  $A,B\in \mathcal H^{\tilde P}$.
Since $H'$ and $H''$ are orthogonal we may assume that $A,B\in H'$ or $A,B\in H''$.

In the former case $\tilde \pi_\ast A =\tilde \pi_\ast B =0$ and $K^fA = K A$, $K^f B= K B$, for $A$ and $B$ are tangent to the fibre at $\tilde x$.
Thus, $\tilde{\bf h}(A,B) = \tilde{\bf h}_f(A,B).$

In the latest case, $K^fA = 0$, $K^f B= 0$, and  $\tilde \pi_\ast A$ and $\tilde \pi_\ast B$ are members of $\mathcal H^{P}$.
Thus $\tilde{\bf h}_f(A,B) =\tg(\pi_\ast A, \pi_\ast B) = \tilde{\bf h}(A,B)$.
Consequently, $\tilde{\bf h}_f = \tilde{\bf h}$ on $\mathcal H^{\tilde P} \times \mathcal H^{\tilde P}$.

If $A,B\in V_{\tilde \xi}$ then  $\tilde \pi_\ast A$ and $\tilde \pi_\ast B$ are tangent to the fibre $P^{-1}(x)$ at $\tilde x$.  We have
\begin{align*}
\tilde{\bf h}_f(A,B) & =\tg_f(\tilde \pi_\ast A,\tilde \pi_\ast B) = f^2(x) \tg ( \tilde \pi_\ast A,\tilde \pi_\ast B) \\
 &=( f\circ \pi)^2 (\xi) \tg ( \tilde \pi_\ast A,\tilde \pi_\ast B) =
 \hat f^2 (\xi) \tilde{\bf h}(A,B).
\end{align*}
Thus $\tilde{\bf h}_f = \hat f^2 \tilde{\bf h}$
on $\mathcal V^{\tilde P} \times  \mathcal V^{\tilde P}$.
Since $\mathcal V^{\tilde P}$ and $\mathcal H^{\tilde P}$ are mutually orthogonal,
$\tilde{\bf h}_f = \hat f^2 \tilde{\bf h}$ on  $\mathcal V^{\tilde P} \times  T(\tilde E^1)$.
\end{proof}

After these preparations we can state  a collapse theorem for a unit horizontal bundle.

\begin{thm}\label{Thm_Col_Thm}  Let $P:(\tilde M, \tilde g) \to (M, g)$ be a Riemannian submersion with $\tilde M$ compact and integrable horizontal distribution. Suppose that $(f_n: M \to (0,+\infty))_{n \in \mathbb{N}}$
is an uniformly bounded sequence of warping functions. We equip each unit  horizontal bundle $\tilde E^1_{f_n}$ and the sphere bundle $SM$ with the
$(p,q)$-metric. Let $\tilde M_n = (\tilde M, \tilde g_{f_n})$ and $\tilde E^1_n = (\tilde E^1_{f_n},\tilde{\bf h}_{f_n}) $.
Then  $\lim_{n\to \infty} \tilde E^1_{n} = SM$ iff
$\lim_{n \to \infty} \tilde M_{n}=M$.
\end{thm}

\begin{proof} From Claim \ref{CTClm3} it follows that $\tilde P : (\tilde E^1_{f_n}, \tilde{\bf h}_{f_n})\to (SM, {\bf h})$
is a Riemannian submersion.
We put $\hat f_n = f_n \circ \pi$, where $\pi:SM\to M$ is the natural projection.

Suppose that $\lim \tilde M_{n} = M $.
Take $\varepsilon >0$. By Theorem  \ref{Thm_Szymon}
($\Rightarrow$) there exist $N>0$ such that for
every $n>N$ there exist an $(\varepsilon/2)$-net  $A^{(n)} = \{x_1,\dots, x_k\} \subset M$ and $f_n| A^{(n)} < \varepsilon/2$.
For each $i=1,\dots, k$ there exist $\xi_{i j}$, $j=1,\dots, l$ such that $\{\xi_{i1},\dots , \xi_{i l}\}$ is an $(\varepsilon/2)$-net in the fibre
$S_{x_i}$.  Put $\hat A^{(n)} = \{\xi_{ij}: i=1,\dots, k;\,\, j=1,\dots l\}$.  By Lemma \ref{net_i_ii}(i), $\hat A^{(n)} $ is an $\varepsilon$-net in $SM$.
Since $\hat f_n(\xi_{ij}) =f_n(x_i)<\varepsilon$, $\hat f_n | \hat A^{(n)} <\varepsilon.$ Consequently. $\lim \tilde E^1_{n} = SM$ by
Theorem \ref{Thm_Szymon} ($\Leftarrow$).

Conversely, suppose that $\lim \tilde E^1_{n} = SM$
Let $\varepsilon >0$. By Theorem \ref{Thm_Szymon} ($\Rightarrow$) there exists $N>0$ such that for every $n>N$ there exists
an $\varepsilon$-net $ \hat A^{(n)} = \{ \zeta_1,\dots, \zeta_m\}$ such that $\hat f_n |\hat A^{(n)} < \varepsilon$.
Then, by Lemma \ref{net_i_ii}(ii) $\pi\big(A^{(n)} \big)$
is an $\varepsilon$-net in $M$. Moreover, for every
$j=1,\dots,m$, $f(\pi(\zeta_j))=\hat f(\zeta_j)<\varepsilon$. Thus $f|A^{(n)} <\varepsilon.$ Consequently, $\lim \tilde M_{n}=  M$, by Theorem \ref{Thm_Szymon}
($\Leftarrow$).
\end{proof}

\vskip20pt
\noindent{\sc Wojciech Koz\l owski}\\
Faculty of Mathematics and Computer Science, University of \L\' od\' z\\
ul. Banacha 22, 90-238 \L\' od\' z, Poland\\
e-mail: {\tt wojciech@math.uni.lodz.pl}
\vskip10pt
\noindent{\sc Szymon M. Walczak}\\
Faculty of Mathematics and Computer Science, University of \L\' od\' z\\
ul. Banacha 22, 90-238 \L\' od\' z, Poland\\
e-mail: {\tt sajmonw@math.uni.lodz.pl}
\vskip20pt
\noindent {\it 2000 Mathematics Subject Classification}: 53B21, 70G45.

\noindent {\it Key words and phrases}:  Cheeger-Gromoll metric, Gromov-Hausdorff distance, Riemannian submersion, unit tangent bundle.
\end{document}